\documentclass[aos,MSNbibl,seceqn,dvips]{arximspdf}
\usepackage{dcolumn}
\usepackage{graphicx}

\doi{10.1214/13-AOS1165} 
\volume{41}
\issue{5}
\pubyear{2013}
\firstpage{2639}
\lastpage{2667}

\makeatletter
\newcolumntype{d}[1]{D{.}{.}{#1}}
\newcommand{\rrVert}{\Vert}
\newcommand{\llVert}{\Vert}
\newtheorem{theorem}{Theorem}

\newtheorem{proposition}{Proposition}
\newcommand{\bX}{\mathbf{X}}
\newcommand{\Xbf}{\mathbf{X}}
\newcommand{\bY}{\mathbf{Y}}
\newcommand{\Zbf}{\mathbf{Z}}
\newcommand{\ubf}{\mathbf{u}}
\makeatother

\begin{document}
\begin{frontmatter}

\title{High-dimensional influence measure}
\runtitle{High-dimensional influence measure}\vspace*{6pt}

\begin{aug}
\author[A]{\fnms{Junlong} \snm{Zhao}\thanksref{t1}\ead[label=e1]{zhaojunlong928@126.com}},
\author[B]{\fnms{Chenlei} \snm{Leng}\corref{}\thanksref{t2}\ead[label=e2]{C.Leng@warwick.ac.uk}},
\author[C]{\fnms{Lexin} \snm{Li}\thanksref{t3}\ead[label=e3]{li@stat.ncsu.edu}}
\and\break
\author[D]{\fnms{Hansheng}~\snm{Wang}\thanksref{t4}\ead[label=e4]{hansheng@gsm.pku.edu.cn}}\vspace*{6pt}
\runauthor{Zhao, Leng, Li and Wang}
\affiliation{Beihang University,
University of Warwick and National University of Singapore,
North Carolina State University and Peking University}
\address[A]{J. Zhao\\
School of Mathematics\\
\quad and System Science\\
Beihang University\\
Beijing\\ 
China\\
\printead{e1}}
\address[B]{C. Leng\\
Department of Statistics\\
University of Warwick\\
Coventry\\ 
United Kingdom\\
and\\
Department of Statistics\\
\quad and Applied Probability\\
National University of Singapore\\
\printead{e2}\vspace*{6pt}}
\address[C]{L. Li\\
Department of Statistics\\
North Carolina State University\\
Raleigh, North Carolina\\
USA
\printead{e3}}
\address[D]{H. Wang\\
Guanghua School of Management\\
Peking University\\
Beijing\\
China\\
\printead{e4}}
\end{aug}\vspace*{6pt}
\thankstext{t1}{Supported by National Natural Science Foundation of
China (11101022)
and Ministry of Education Humanities and Social Science Foundation
Youth project (10YJC910013).}
\thankstext{t2}{Supported by National University of Singapore research grants.}
\thankstext{t3}{Supported by NSF Grant DMS-11-06668.}
\thankstext{t4}{Supported by National Natural Science Foundation of
China (11131002, 11271032),
Fox Ying Tong Education Foundation,
the Business Intelligence Research Center of Peking University
and the Center for Statistical Science of Peking University.}

\received{\smonth{7} \syear{2013}}

%
\begin{abstract}
Influence diagnosis is important since presence of influential
observations could lead to distorted analysis and misleading
interpretations. For high-dimensional data, it is particularly so, as
the increased dimensionality and complexity may amplify both the chance
of an observation being influential, and its potential impact on the
analysis. In this article, we propose a novel high-dimensional
influence measure for regressions with the number of predictors far
exceeding the sample size. Our proposal can be viewed as a high-dimensional
counterpart to the classical Cook's distance. However,
whereas the Cook's distance quantifies the individual observation's
influence on the least squares regression coefficient estimate, our new
diagnosis measure captures the influence on the marginal correlations,
which in turn exerts serious influence on downstream analysis including
coefficient estimation, variable selection and screening. Moreover, we
establish the asymptotic distribution of the proposed influence measure
by letting the predictor dimension go to infinity. Availability of this
asymptotic distribution leads to a principled rule to determine the
critical value for influential observation detection. Both simulations
and real data analysis demonstrate usefulness of the new influence
diagnosis measure.
\end{abstract}

%

\end{frontmatter}

\section{Introduction}\label{sec1}

An observation is flagged influential if some important features of the
analysis are substantially altered after this observation is removed
\cite{Cook:1979}. Presence of influential observations would possibly
lead to distorted analysis and misleading results
\cite{Draper:Smith:1998}, and therefore it is important to be alert to
influential observations and take them into consideration when
interpreting the results. In the classical normal linear model setup,
regression coefficient estimate was chosen, naturally, as the feature
whose substantial change defines influential observations. Toward that
end, \cite{Cook:1977} proposed a difference measure between the OLS
estimate on the full data and that on the subset of data without the
observation in question. This measure, which is later on referred in
the statistical literature as the \emph{Cook's distance}, quantifies
the contribution, or influence, of individual data observation on the
regression coefficient estimate. Consequently an observation with a
large Cook's distance is deemed as influential. Since its introduction,
the Cook's distance has been routinely employed in regression analysis,
due to its clear interpretation from the case deletion point of view,
and its easy computation without having to re-estimate the model for
each removed observation. The topic is covered in most standard
regression textbooks, and it is implemented in popular statistical
software such as R and SAS.

The problem of influence diagnosis has since attracted considerable
attention and been systematically investigated for various models and
analyses. Examples include linear regression models
\cite{Cook:1977,Chatterjee:Hadi:1988,Cook:Weisberg:1982}, categorical
data analyses \cite{Anderson:1992}, generalized linear models
\cite{Williams:1987,Thomas:Cook:1989,Davison:Tsai:1992}, generalized
estimation equations \cite{Preisser:Qaqish:1996}, linear mixed models
\cite{Christensen:Pearson:Johnson:1992,Banerjee:Frees:1997,Banerjee:1998},
generalized linear mixed models \cite{Xiang:Tse:Lee:2002},
semiparametric mixed models \cite{Fung:Zhu:Wei:He:2002}, growth curve
models \cite{Pan:Fang:2002}, incomplete data analysis
\cite{Zhu:Lee:Wei:Zhou:2001}, perturbation theory
\cite{Chritchley:Atkinson:Lu:Biazi:2001,Zhu:Ibrahim:Lee:Zhang:2007,Zhu:Ibrahim:Cho:2012},
among others. For an excellent review on the latest developments in the
field of influence diagnosis, we refer to \cite{Zhu:Ibrahim:Cho:2012}.

Thanks to the aforementioned works, substantial insights have been
gained on influence diagnosis. However, it is important to note that,
all existing diagnosis approaches have been developed under the
assumption that the number of predictors in regression is fixed. As
such, none is immediately applicable to high-dimensional regression
analysis, where the number of predictors $p$ far exceeds the sample
size $n$. On the other hand, nowadays prevailing in both science and
business are data with unprecedented size and dimensionality, calling
for the development of high-dimensional influence diagnosis. Detection
of influential observations in high-dimensional data analysis, in our
opinion, is equally, or to some extent, even more important than in a
classical setup. This is partly because the increased dimensionality
and complexity of the data may amplify both the chance of an
\mbox{observation} being influential as well as its potential impact on
the analysis. Moreover, the peculiar data observations themselves may
be of practical importance in addition to data modeling. The diagnosis
task, nevertheless, is more challenging in high-dimensional data
analysis, and is far from a direct extension of existing diagnosis
approaches. To the best of our knowledge, influence diagnosis in a
high-dimensional setting has received little attention despite its
evident importance.

The first challenge is the definition of influential observation. In
other words, which feature of the analysis should one choose such that
its substantial alternation defines an influential observation? In the
classical setup, an observation is deemed influential if it incurs
serious change in regression coefficient estimate. In high-dimensional
regression where $p > n$, the ordinary least squares estimator is
highly unstable as the gram matrix is not invertible. On the other
hand, we recognize that variable selection and variable screening are
of particular importance in high-dimensional regression analysis. There
has been a vast literature on variable selection in recent years,
including the LASSO \cite{Tibshirani:1996}, the adaptive LASSO
\mbox{\cite{Zou:2006,Zhang:Lu:2007,Wang:Leng:2007}}, the SCAD
\cite{Fan:Li:2001}, the bridge estimator
\cite{Fu:1998,Huang:Horowitz:Ma:2007}, the LARS
algorithm~\cite{Efron:Hastie:Johnstone:Tibshirani:2004}, the Dantzig
selector \cite{Candes:Tao:2007}, the sure independence screening rule
\cite{Fan:Lv:2008}, SIS, the forward regression \cite{Wang:2009}, FR,
among many others. Underlying all those selection methods, one
statistic plays a critical role and, that is, the \emph{marginal}
covariance, or equivalently, \emph{marginal} correlation between the
response and the individual covariates. To clarify, we note that, SIS
is directly defined based on this statistic, whereas the first step of
the forward regression hinges on the estimated marginal covariance too.
In addition, the sample marginal covariance, in addition to the Gram
matrix, is an important input for the well celebrated LARS algorithm,
as well as the LASSO, the adaptive LASSO and the Dantzig selector.

Motivated by this vital observation, we choose the marginal correlation
as the feature that defines influential observation. We propose a new
influence diagnosis measure, which continues to utilize the
leave-one-out idea of the classical Cook's distance, but is based on
the combined marginal correlations between the response and all
predictors. The new measure is applicable to high-dimensional setting
where $p > n$, and is very fast and easy to compute. Unlike the
classical Cook's distance that quantifies the individual observation's
influence on the least squares coefficient estimate, the new measure
captures the influence on the marginal correlation, which \emph{in
turn} exerts serious impact on variable selection and other downstream
analysis. The choice of the marginal correlation as the defining
feature of our influence diagnosis does not imply that the marginal
correlation is our ultimate goal of interest. Instead, it reflects
influence on important analysis features including parameter
estimation, variable selection and screening. This definition of
influential observation in a high-dimensional setting can be viewed as
our first contribution.

Our second contribution is that the explicit asymptotic distribution
for the proposed influence measure is derived. Availability of this
asymptotic theory offers a principled guidance to determine the
critical value for the influence measure. Subsequently, we propose a
false discovery rate based procedure for that purpose
\cite{Benjamini:Hochberg:1995,Benjamini:Hochberg:2000}. We remark that,
in the classical setup where $p$ is fixed, a standard Taylor's
expansion type analysis \cite{Cook:1977} revealed that the classical
Cook's distance's major variability is due to the observation under
investigation and its sample size is only one. This rules out the
possibility of establishing a standard asymptotic theory for the
classical Cook's distance. To determine an appropriate threshold value
for the classical Cook's distance, its distribution can be obtained by
bootstrap if the true model is a parametric linear model. However, such
a bootstrap procedure requires a parametric model assumption and can be
computationally expensive especially for high-dimensional data. By
contrast, the asymptotic distribution of the proposed influence measure
is attainable in our setup, since the predictor dimension goes to
infinity along with the sample size, and the threshold is easy to
obtain.

When facing high-dimensional data diagnosis, an intuitive solution is
to continue using the classical Cook's distance but to replace the OLS
coefficient estimate with a regularized estimate, for instance, a LASSO
estimate. This modified Cook's distance approach could be particularly
useful when data perturbation concentrates on the nonzero coefficients,
as it avoids unnecessary variability caused by irrelevant covariates.
However, it also has several limitations. First, this solution
interweaves influence diagnosis with variable selection, which can be
flawed if the influence is reflected on variable selection itself. For
instance, an influential observation may substantially alter the chosen
tuning parameter of the LASSO, resulting in a totally different
regularized coefficient estimate, which in turn affects the modified
Cook's distance. Second, the tuning parameter of the LASSO, in
principle, should be updated for every reduced data set, and this
re-estimation requirement can be very expensive computationally,
especially when the regression dimension $p$ is large. Third, the
asymptotic properties of the modified Cook's distance seem intractable
analytically, which makes the thresholding of influential data
difficult, whereas a bootstrap alternative to choose the thresholding
value is again computationally expensive. Moreover, while there exist
many competing variable selection methods, it is unclear which
selection method is the best choice in the context of influence
diagnosis. By contrast, our influence measure is not constrained by any
particular variable selection method, and this flexibility could
benefit downstream analysis. In Section~\ref{sec3}, we carry out an
intensive numerical study to compare this modified Cook's distance with
our proposal, and this detailed comparison can be viewed as the third
contribution of this article.

\begin{figure}

\includegraphics{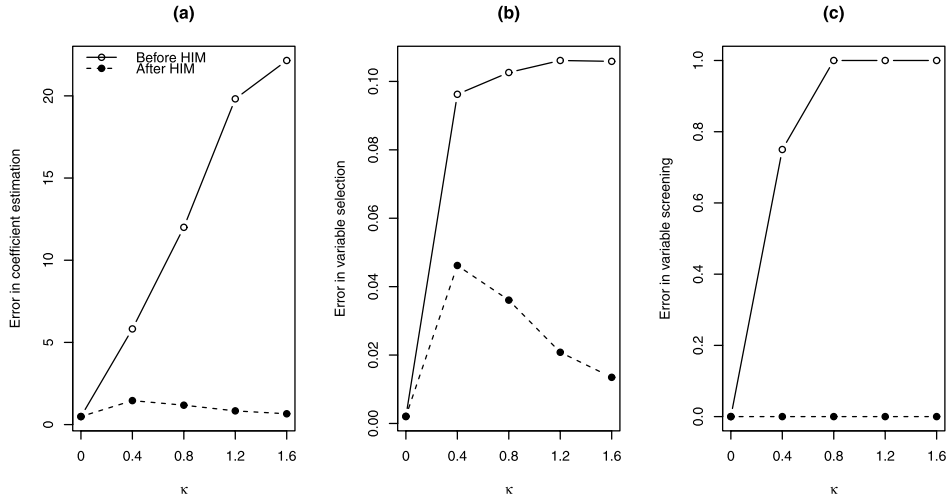}

\caption{Effect of influential points on parameter
estimation \textup{(a)}, variable selection \textup{(b)} and variable screening \textup{(c)}, as
the perturbation parameter $\kappa$ varies. ``Before HIM'' denotes the
analysis on the full data, and ``After HIM'' denotes the analysis on the
reduced data after removing the influential observations flagged by our
proposed high-dimensional measure (HIM).}\label{effect}
\end{figure}

Before we proceed, we quickly show a simulated example to illustrate
two points: first, how various aspects of a high-dimensional regression
analysis, \mbox{including} regression coefficient estimation, variable
selection and variable screening, can be seriously affected by
influential observations, and second, how our proposed measure can help
limit such influence. The data was generated from a linear model with
$p=1000$ predictors, $n=100$ observations, among which 10 observations
were influential. The magnitude of the influence was dictated by a
scalar $\kappa$ with a larger value indicating a larger influence. More
details can be found in the setup of model~1 in Section~\ref{sec3}.
Evaluations include error in coefficient estimation, error in variable
selection after applying the LASSO \cite{Tibshirani:1996}, and error in
variable screening after applying the SIS \cite{Fan:Lv:2008}. The
results are averaged over 200 simulation replicates, and are reported
in Figure~\ref{effect}. It is clearly seen from the plot that,
influential observations could have drastic effects on various features
for high-dimensional data analysis. Meanwhile, our marginal correlation
based diagnosis could greatly help control the adverse effects after
detecting and removing those influential data points.

%
%

The rest of the artlicle is organized as follows. Section~\ref{sec2}
begins with a review of the classical Cook's distance, then presents
our new high-dimensional influence measure, along with a comparison
with the Cook's distance, the asymptotic properties and a power study.
Section~\ref{sec3} includes an intensive simulation study and a
microarray data analysis. Section~\ref{sec4} presents a generalization
of our proposal from the normal linear model to the generalized linear
model. Section~\ref{sec5} concludes the paper with a discussion. All
technical proofs are given in the \hyperref[sec6]{Appendix} and the supplementary
material \cite{suppJZ}.

\section{High-dimensional influence measure}\label{sec2}

\subsection{Linear models and classical Cook's distance}\label{sec2.1}

In this article, we focus on influence diagnosis in the context of the
classical linear regression model. Meanwhile, we note that the proposed
idea can be readily extended to a much broader class of regression
models, and we will discuss one such extension in Section~\ref{sec4}.
Consider the following model:
%
%
\begin{equation}
\label{eq:model} Y_i= \beta_0 + \Xbf_i^\top
\bolds{\beta}_1+\varepsilon_i,
\end{equation}
where the pair $(Y_i,\Xbf_i)$, $1\leq i\leq n$, denote the observation
of the $i$th subject, $Y_i\in\mathbb{R}$ is the response variable,
$\Xbf_i=(X_{i1},\ldots,X_{ip})^\top\in\mathbb{R}^p$ is the associated
\mbox{$p$-}dimensional predictor vector, and $\varepsilon_i\in\mathbb{R}$ is a
mean zero normally distributed random noise. Let $\bolds{\beta}=
(\beta_0, \bolds{\beta}_1^\top)^\top$ denote the coefficient vector.
Under the classical setup of $n > p$, the OLS estimate of
$\bolds{\beta}$ is obtained by minimizing the objective function
$\sum_{i=1}^n (Y_i - \beta_0 - \Xbf_i^\top\bolds{\beta}_1)^2$, and the
solution is $\hat{\bolds{\beta}} = (\bX^\top\bX)^{-1}\bX^\top\bY$,
where $\bY=(Y_1,\ldots, Y_n)^\top$ denotes the $n \times1$ response
vector, and $\bX$ denotes the $n \times(p + 1)$ design matrix with the
$i$th row being $p+1$ dimensional vector $(1, \Xbf_i^\top)$,
$i=1,\ldots,n$.

To quantify the influence of the $k$th observation on regression,
$1\leq k \leq n$, \cite{Cook:1977} employed the leave-one-out idea by
studying the OLS estimate of $\beta$ while the $k$th observation is
excluded from estimation. That is, one minimizes the modified objective
function $\sum_{i=1,i\neq k}^{n}
(Y_i-\beta_0-\Xbf_i^\top\bolds{\beta}_1)^2$. The new estimate is of the
form $\hat{\bolds{\beta}}{}^{(k)} = (\bX_{(k)}^\top\bX_{(k)})^{-1}
\bX_{(k)}^\top\bY_{(k)}$, where $\bY_{(k)}$ is the $(n-1) \times1$
response\vspace*{1pt} vector with $Y_k$ removed, and $\bX_{(k)}$ is the
$(n-1)\times(p+1)$ design matrix with the $k$th row $\Xbf_k$ removed.
Cook \cite{Cook:1977} naturally chose the estimate of $\bolds{\beta}$
to define influence, and intuitively, if an observation is influential,
the difference between $\hat{\bolds{\beta}}$ and
$\hat{\bolds{\beta}}{}^{(k)}$ is expected to be large. This leads to
the following discrepancy measure, that is, the Cook's distance:
%
%
%
\begin{equation}
\label{eq:CD} D_k=\frac{ \{
\hat{\bolds{\beta}}{}^{(k)}-\hat{\bolds{\beta}} \}^\top\bX^\top\bX\{
\hat{\bolds{\beta}}{}^{(k)}-\hat{\bolds{\beta}} \}} {(p+1)\hat\sigma{}^2},
\end{equation}
where $\hat\sigma{}^2=(n-p-1)^{-1}\sum_{i=1}^{n}
(Y_i-\hat\beta_0-\Xbf_i^\top\hat{\bolds{\beta}})^2$.

In the high-dimensional regression setting, the classical Cook's
distance (\ref{eq:CD}) encounters some difficulties. When $p$ is close
to $n$, the OLS estimate is known to be unstable, which would in turn
cause $D_k$ to be unstable. When $p > n$, the classical Cook's distance
is not directly computable, because the OLS estimator $\hat\beta$
becomes unstable. For those reasons, the regression coefficient
estimate may no longer be the best choice to define influence in
high-dimensional analysis. This motivates us to consider an alternative
influence measure for high-dimensional data.

\subsection{High-dimensional influence measure}\label{sec2.2}

In high-dimensional regression analysis where $p \approx n$ or $p > n$,
variable selection (screening) plays a central role, whereas marginal
covariance or correlation is crucial to the majority of variable
selection approaches. Motivated by this observation, for
high-dimensional data influence diagnosis, we choose marginal
correlation, instead of regression coefficient, as the feature that
defines influence. Individual observation's influence on marginal
correlation is to transmit to various features of downstream analysis,
such as variable selection and coefficient estimation.

More specifically, we first define the marginal correlation as $\rho_j
= E\{(X_j - \mu_{xj})(Y - \mu_y)\} / (\sigma_{xj} \sigma_y)$, where
$\mu_{xj} = E(X_j)$, $\mu_y = E(Y)$, $\sigma_{xj}^2 =
\operatorname{var}(X_j)$ and $\sigma_y^2 = \operatorname{var}(Y)$. We
then obtain the sample estimate, $\hat{\rho}_j = \{\sum_{i=1}^{n}
(X_{ij} - \hat{\mu}_{xj})(Y_i - \hat{\mu}_y)\}/\{n \hat{\sigma}_{xj}
\hat{\sigma}_y \}$, for $j = 1,\ldots,p$, where $\hat{\mu}_{xj},
\hat{\mu}_y, \hat{\sigma}_{xj}$ and $\hat{\sigma}_y$ are the sample
estimates of $\mu_{xj}$, $\mu_y$, $\sigma_{xj}$ and $\sigma_y$,
respectively. Next, we continue to use the leave-one-out principle as
in the classical Cook's distance case, and compute the marginal
correlation with the $k$th observation removed as
\[
\hat{\rho}{}_j^{(k)} = \frac{\sum_{i=1, i \neq
k}^{n} (X_{ij} - \hat{\mu}{}^{(k)}_{xj} ) (Y_i -
\hat{\mu}{}^{(k)}_y )}{(n-1) \hat{\sigma}{}^{(k)}_{xj}
\hat{\sigma}{}^{(k)}_y },\qquad j =1, \ldots, p, k = 1, \ldots, n,
\]
where $\hat{\mu}{}_{xj}^{(k)}, \hat{\mu}{}_y^{(k)},
\hat{\sigma}{}_{xj}^{(k)}$ and $\hat{\sigma}{}_y^{(k)}$ are the
corresponding sample estimates with the $k$th observation removed.
Finally, we define the influence measure based on the marginal
correlation as
%
%
\begin{equation}\label{D0DK}
\mathcal{D}_k=\frac{1}{p} \sum
_{j=1}^p \bigl(\hat\rho_j-\hat\rho{}_{j}^{(k)} \bigr)^2.
\end{equation}
We refer to $\mathcal{D}_k$ as the \emph{high-dimensional influence
measure}, or HIM for brevity. We make a few remarks. First, we note
that the marginal correlation can be easily computed regardless of the
predictor dimension, and such computational advantage is practically
very useful for high-dimensional data analysis. Second, the proposed
influence measure is built upon the marginal correlation coefficient,
and is effectively scale invariant. However, it does \emph{not} imply
that marginal correlation is the ultimate feature of interest in our
influence diagnosis. Instead, a substantial change on the marginal
correlation caused by a data point is to exert influence on important
features such as variable selection and parameter estimation, as we
have seen in Figure~\ref{effect}. As such, for an estimation method to
be robust to unexpected \emph{perturbation}
\cite{Chritchley:Atkinson:Lu:Biazi:2001,Zhu:Ibrahim:Lee:Zhang:2007,Zhu:Ibrahim:Cho:2012},
the sample marginal correlation should be sufficiently robust. This is
an important and necessary condition, although not necessarily
sufficient. Finally, use of the \emph{marginal} correlation to define
the influence measure does \emph{not} imply that we assume a
\emph{marginal model}. Instead, we still assume the joint
model~(\ref{eq:model}). As it may seem unclear how a marginal measure can
capture the influence for a joint model, we will demonstrate through a
simple joint model later in Section~\ref{sec2.5} that, the newly
defined $\mathcal{D}_k$ can indeed identify the influential observation
with probability one. This use of marginal correlation is also similar
in spirit to the sure independence screening procedure for a joint
normal model \cite{Fan:Lv:2008}, but is in a different context. Fan and
Lv \cite{Fan:Lv:2008} use marginal correlation for the variable
screening purpose, while we use it for influence diagnosis.

The proposed high-dimensional influence measure also shares some
similarity as the classical Cook's distance. Note that the Cook's
distance can be reformulated~as
%
\begin{equation}\label{cookdis1}
D_k=\frac{\hat{\epsilon}_k^2}{p
\hat{\sigma}^2}\frac{h_{kk}}{(1-h_{kk})^2},\qquad k=1,\ldots,n,
\end{equation}
where $\hat\epsilon_k=\hat Y_k-Y_k$ is the residual and
$h_{kk}=\Xbf_k^\top(\bX^\top\bX)^{-1}\Xbf_k$, $k=1,\ldots, n$ is the
$(k)$th diagonal element of the hat matrix
$\bX(\bX^\top\bX)^{-1}\bX^\top$. Clearly, $D_k$ is an increasing
function of both $|\hat\epsilon_k|$ and $h_{kk}$. As such, an
observation has a large value in Cook's distance, if it has a large
residual or it is a high leverage point in terms of~$h_{kk}$. Our
proposed information measure shares a similar spirit. In
Section~\ref{sec2.3}, we will derive a decomposition of our influence
measure $\mathcal{D}_k$ under some conditions, and will show that
$\mathcal{D}_k$ is mainly dominated by a term called $B_2$, which is of
the form
\[
B_2=\frac{(n-2)}{pn(n-1)^2}\sum_{j=1}^p
Y_k^2X^2_{kj} =\frac{(n-2)}{pn(n-1)^2}
Y_k^2 \|\Xbf_k\|^2.
\]
Consequently the $k$th data point $(\Xbf_k, Y_k)$ is more likely to
be marked influential, if it has a large response and a large value of
$\|\Xbf_k\|^2$. Here $\|\Xbf_k\|^2$ plays a similar role as $h_{kk}$ in
the classical Cook's distance, for detecting influential points induced
mainly by covariates, whereas $Y_k$ plays a similar role as the
residual in the Cook's distance, for detecting the influential point
induced by abnormal responses.

\subsection{Theoretical properties}\label{sec2.3}

We next establish the asymptotic distribution of the proposed high-dimensional influence
measure $\mathcal{D}_k$ as both the sample size
$n$ and
the dimensionality $p$ go to infinity. Toward that end, we impose the
following conditions.
\begin{longlist}[(C.3)]
\item[(C.1)] For any fixed $j = 1, \ldots, p$, $\rho_j$ is constant
    and does not change as $p$ increases.

\item[(C.2)] For the covariance matrix $\bolds{\Sigma}= \operatorname
    {cov}(\Xbf)$, with the eigen-decomposition
    $\bolds{\Sigma}=\sum_{j=1}^p \lambda_j \ubf_j \ubf_j^\top$, it is
    assumed that $l_p=\sum_{j=1}^{p} \lambda_j^2 = O(p^r)$ for some
    \mbox{$0\leq r<2$}.

\item[(C.3)] The predictor $X_i$ follows a multivariate normal
    distribution and the random noise $\varepsilon_i$ follows a normal
    distribution.
\end{longlist}
Condition \textup{(C.1)} is very general, since it only requires that
for any fixed $j$, $\rho_j$ is a constant independent of $p$. A
sufficient condition for condition~\textup{(C.2)} to hold is that all
eigenvalues of $\bolds{\Sigma}$ are finite. This condition also permits
eigenvalues of $\bolds{\Sigma}$ to diverge to infinity but at a slower
rate compared to the dimensionality. The normality assumption on $\Xbf$
is mainly for convenience, and can be relaxed, for instance, to
distributions with sub-Gaussian tails, at the expense of more lengthy
proofs. In addition, since the error term is assumed normal, $Y$ is
normally distributed.

Next, we derive a decomposition of $\mathcal{D}_k$, that is, to serve
as a basis for its asymptotic distribution. The result is presented in
a way such that $\mu_y$, $\mu_{xj}$ are assumed to be 0 and $\sigma_{xj}$,
$\sigma_y$ are 1 for $1\leq j\leq p$. This leads to simplified estimates
$\hat\rho_j=n^{-1}\sum_{1\leq i\leq n} X_{ij}Y_i$ and
$\hat\rho{}_j^{(k)}=n^{-1}\sum_{i\not=k} X_{ij}Y_i$. On the other hand,
we note that this standardization is only for the purpose of
simplifying the presentation and it loses no generality. As we will
show later in Proposition~\ref{Dksample}, \mbox{replacing} the unknown
quantities $\mu_{xj}$, $\mu_y$, $\sigma_{xj}$ and $\sigma_y$ with
their consistent sample estimates would not alter $\mathcal{D}_k$'s
asymptotic distribution. For $t,s=1,\ldots,n$, let $K_{p,ts}=\sum_j
X_{tj}X_{sj}/p$ and $c_p = \max_{1 \leq j \leq p} \lambda_j$. After
some algebraic computation, we obtain that
%
%
%
\begin{eqnarray}\label{D_k_decom}
\qquad \mathcal{D}_k&=&
\frac{1}{p}\sum_{j=1}^p \Biggl\{
\frac{1}{n(n-1)}\sum_{1\leq t\leq n}^{t\neq k}
Y_tX_{tj}-\frac{1}{n} Y_kX_{kj}
\Biggr\}^2\nonumber
\\
&=&\frac{1}{\{n(n-1)\}^2}\sum_{t=1}^n
Y_t^2 K_{p,tt}+\frac{(n-2)}{n(n-1)^2}Y_k^2K_{p,kk}
\nonumber\\[-8pt]\\[-8pt]
&&{} +\frac{1}{[n(n-1)]^2} \sum_{t\neq s} Y_tY_s K_{p,ts}
-\frac{2}{n(n-1)^2} \sum_{t=1,t\neq k}^n
Y_kY_t K_{p,tk}\nonumber
\\
&:=& B_1+B_2+B_3-2B_4.\nonumber
\end{eqnarray}
Then we have the following result on the expectation of $\mathcal{D}_k$ along
with the variance of its decomposition in terms of $B$'s.

%
\begin{proposition}\label{i-cook}
Suppose that $(\Xbf_i,Y_i)$ are i.i.d. observations and that
\textup{(C.1)}~and~\textup{(C.3)} hold. Then it holds that
\[
E(\mathcal{D}_k)=\bigl[n(n-1)\bigr]^{-1}E
\bigl(Y_k^2\bigr)E(K_{p,kk})+O
\bigl(n^{-2}p^{-1}l_p^{1/2}\bigr).
\]
In addition, $\operatorname{var}(B_1)=O(n^{-7})$,
$\operatorname{var}(B_2)=O(n^{-4})$, $\operatorname{var}(B_3)=O(c_p^2
n^{-5} p^{-2})+O(p^{-2}n^{-6})$ and $\operatorname{var}(B_4)=O(l_p
p^{-2}n^{-5})+O(c_p^2p^{-2}n^{-4})$.
\end{proposition}

Now we return to the asymptotic distribution of $\mathcal{D}_k$.
Proposition~\ref{i-cook} helps to derive the asymptotic distribution of
$\mathcal{D}_k$. We first present the result assuming $\mu_{xj}$,
$\mu_y$, $\sigma_{xj}$ and $\sigma_y$ are all known. Then we obtain the
asymptotic distribution when $\mu_{xj}$, $\mu_y$, $\sigma_{xj}$ and
$\sigma_y$ are replaced by their sample estimates.

%
\begin{theorem}\label{Dk} Suppose that \textup{(C.1)--(C.3)} hold. When there
is no influential point and $\min\{n,p\}\longrightarrow\infty$, we have
\[
n^2\mathcal{D}_k \longrightarrow\chi^2(1),
\]
where $\chi^2(1)$ is the chi-square distribution with one degrees of
freedom.
\end{theorem}

Next, we consider the asymptotic distribution of $\mathcal{D}_k$ when
$\mu_{xj}$, $\mu_y$, $\sigma_j$ and $\sigma_y$ are unknown. A natural
choice is to replace them by their corresponding sample moment
estimates\vspace*{-1pt} as $\hat{\mu}_y=\sum_i Y_i/n$, $\hat{\mu}_{xj}=\sum_i
X_{ij}/n$, $\hat{\sigma}_{xj}^2=\sum_i
(X_{ij}-\hat{\mu}{}_{xj})^2/(n-1)$ and $\hat{\sigma}{}_y^2=\sum_i (Y_i-
\hat{\mu}_y)^2/(n-1)$. Another choice is to employ robust estimators,
for example, the median in place of the mean, and the median absolute
deviation in place of the standard deviation. The following proposition
shows that the conclusion of Theorem~\ref{Dk} continues to hold as long
as $u_{xj}$, $u_y$, $\sigma_{xj}$ and $\sigma_y$ are replaced by
$\sqrt{n}$-consistent estimates  under certain moment assumptions. Let
$\dot{Y}_t=(Y_t-\mu_y)/\sigma_y$,
$\dot{X}_{tj}=(X_{tj}-u_{tj})/\sigma_{tj}$, $t=1,\ldots,n,j=1,\ldots,
p$ and $(Q_{xj},R_{xj})=((\hat{\mu}_{xj}-\mu_{xj})/\sigma_{xj},
\sigma_{xj}/\hat{\sigma}_{xj})$ and $(Q_{y},R_{y})$ are defined
similarly. Furthermore,\vspace*{2pt} let $S_{Qx}=\limsup_{n\rightarrow\infty}
E(n^{1/2}Q_{x1})^8$, $S_{Rx}=\limsup_{n\rightarrow\infty}
E[n^{1/2}(R_{x1}-1)]^8$,\vspace*{2pt} $S_{Qy}=\limsup_{n\rightarrow\infty}
E(n^{1/2}Q_{y})^8$ and $S_{Ry}=\break\limsup_{n\rightarrow\infty}
E[n^{1/2}(R_{y}-1)]^8$. We make the following additional assumption.
\begin{longlist}[(C.4)]
\item[(C.4)] For all $1\leq j\leq p$, $(Q_{xj},R_{xj})$ are the
    same symmetric function of $\{\dot{X}_{tj}$, for
    $t=1,\ldots,n\}$; and $(Q_{y},R_{y})$ are also the same
    symmetric function of $\dot{Y}_t$ for $t=1,\ldots,n$. We assume
    that $S_{Qx}$, $S_{Rx}$, $S_{Qy}$ and $S_{Ry}$ are finite.
\end{longlist}

Condition \textup{(C.4)} indicates that, for all $1\leq j\leq p$,
$((\hat{\mu}_{xj}-\mu_{xj})/\sigma_{xj},
\hat{\sigma}_{xj}/\sigma_{xj})$ $=f(\dot{X}_{1j},\ldots,\dot{X}_{nj})$,
where $f(x_1,\ldots,x_p)=(f_1(x_1,\ldots,x_p), f_2(x_1,\ldots,x_p))$
and $f_1$ and $f_2$ are symmetric functions. Condition~\textup{(C.4)}
is a mild condition. Recall that $(\Xbf_i, Y_i)$, $i=1,\ldots, n$ are
i.i.d. normal in Theorem~\ref{Dk}. When
$\hat{\mu}_{xj},\hat{\sigma}_{xj}$ are the moment estimates, we have
$Q_{xj}=n^{-1}\sum_{1\leq t\leq n} \dot X_{tj}\sim N(0,1/n)$ and
consequently $S_{Qx}$ is finite. Moreover, we have $R_{xj}=1/S_{nj}$
where $S_{nj}^2$ is the sample variance of
$\{\dot{X}_{tj},t=1,\ldots,n\}$. Since $S_{n1}^2\sim
\chi^2_{n-1}/(n-1)$, it is easy to verify that $S_{Rx}$ is also finite.
Similarly, $S_{Qy}$ and $S_{Ry}$ are also finite with moment estimates
$\hat\mu_{y}$ and $\hat\sigma_y$. Under the normality of $(\Xbf, Y)$,
\textup{(C.4)} also holds for some robust estimates.

%
\begin{proposition}\label{Dksample}
Assume that $\hat{\mu}_{xj},\hat{\sigma}_{xj}, \hat{\mu}_{y},
\hat{\sigma}_{y}$ are $\sqrt{n}$-consistent and satisfy~\textup{(C.4)}.
Substituting $\mu_{xj}, \mu_y, \sigma_j, \sigma_y$ with their
corresponding estimates in $\mathcal{D}_k$, Theorem~\ref{Dk} continues
to hold under the same conditions.
\end{proposition}
We remark that the asymptotic distribution of the high-dimensional
influence measure $\mathcal{D}_k$ is obtained as the number of
predictor $p$
goes to infinity. This is different from the case of classical Cook's
distance $D_k$ where $p$ is fixed, for which a standard asymptotic
distribution is not attainable. We view this as a \emph{blessing} of
dimensionality in contrast to the usually conceived \emph{curse} of
dimensionality. For more examples of blessing of dimensionality, see
\cite{Donoho:2000} and \cite{Johnstone:2001}.

\subsection{Influence diagnosis}\label{sec2.4}

An important implication of Theorem \ref{Dk} is that we can now obtain
a $p$-value for influence diagnosis. Specifically, for the hypothesis
that the $k$th observation is not influential versus its alternative,
the $p$-value is $P(\chi^2(1) > n^2 \mathcal{D}_k)$. Given that the
number of predictors $p$ is usually large and multiple hypotheses are
tested simultaneously, we employ the false discovery rate based
multiple testing procedure of \cite{Benjamini:Hochberg:1995} to
determine which hypothesis should be rejected while controlling the
family-wise error. Denote $n_{\mathrm{infl}}$ as the number of
influential observations among the $n$ observations, $n_{\mathrm{tp}}$ and
$n_{\mathrm{fp}}$ as the number of the observations that are correctly rejected
and incorrectly rejected, respectively, and $r$ as the total number of
rejections in the $n$ hypotheses testing. Then the power and the false
discovery rate are denoted as $\mathrm{Power} =
n_{\mathrm{tp}}/n_{\mathrm{infl}}$ and $\mathrm{FDR} = n_{\mathrm{fp}} / r$,
respectively. We will set FDR level being small, such as 0.05, and
report the power and other quantities in the numerical study section.
We also remark that more sophisticated alternative multiple testing
procedure, for example, in \cite{Benjamini:Hochberg:2000},
\cite{Efron:Tibshirani:Storey:Tusher:2001} and \cite{Storey:2002}, can
be used in conjunction with our approach, but, that is, not the focus
of this article.

\subsection{A power comparison of two influence measures}\label{sec2.5}

We next study the power property of both the new diagnosis measure and
the Cook's distance via a simple model. This study serves two purposes.
First, we can gain insight about difference between the two diagnosis
measures. Second, it offers evidence that the marginal correlation
based measure is capable of detecting influential observation in a
joint model with a large probability.

More specifically, we consider the model (\ref{eq:model}), but drop the
intercept for simplicity. The predictors $\Xbf_i, i=1,\ldots, n$, are
i.i.d. observations from a multivariate normal distribution
$N_p(0,\bolds{\Sigma})$ where $\bolds{\Sigma}$ is a $p \times p$
covariance matrix with all its diagonal elements $\sigma_{jj}=1$. The
error term $\epsilon_i$ is of the structure $\epsilon_i = e_i + c_i$,
where $e_i$ follows a standard normal distribution and $c_i$ is
constant, $c_2 = \cdots= c_n = 0$. Under this setup, the first
observation is an influential point as long as $c_1 \neq0$, and we aim
to establish the power of both the classical and our proposed
high-dimensional influence measure in identifying this influential
observation. Let $D_i$ be the Cook's distance\vspace*{-1pt} defined in (\ref{eq:CD})
for the $i$th observation, $\mathcal{D}_i^{(c)}$ be the proposed
high-dimensional measure in (\ref{D0DK}), and
$\mathcal{T}_i^{(c)}=n^2\mathcal{D}^{(c)}_i$ be the statistic defined
in Theorem~\ref{Dk}. Moreover, consider the following condition:
\begin{longlist}[(C.5)]
\item[(C.5)] All eigenvalues of $\bolds{\Sigma}$ are positive and
bounded.
\end{longlist}
Then the next theorem states that, both the classical and the
high-dimen\-sional Cook's distance have the power of detecting the
influential observation approaching one under appropriate yet different
conditions.

%
\begin{theorem}\label{Cook-u}
Consider the model stated above.
\begin{enumerate}[2.]
\item Suppose that \textup{(C.1)} and \textup{(C.5)} hold. If
    $\max\{n^{-1}p^6, |c_1|^{-1}n^{2/3}\}\rightarrow0$, then we
    have that for the Cook's distance $D_i$, $P(nD_1-\max_{2\leq
    i\leq n} nD_i>M)\rightarrow1$ for any $M>0$, when $n
    \to\infty$.

\item Suppose\vspace*{-1pt} that \textup{(C.1)} and \textup{(C.2)} hold. If
    $\max\{|c_1|^{-1}(\log n)^{1/2},
    l_pp^{-2}c_1^{-4}n\}\rightarrow0$, then we have that for the
    proposed\vspace*{-1pt} high-dimensional influence measure
    $\mathcal{D}_i^{(c)}$, $P(\mathcal{T}_1^{(c)}-\max_{2\leq i\leq
    n} \mathcal{T}_{i}^{(c)}>M) \to1$ for any $M>0$, when\break \mbox{$\min(n,
    p) \to \infty$}.
\end{enumerate}
\end{theorem}

The proof is given in the supplementary
material \cite{suppJZ}. Here we compare the
two sets of conditions to gain some insight about the difference of the
two diagnosis measures. First,\vspace*{1pt} we examine the condition
$\max\{n^{-1}p^6, |c_1|^{-1}n^{2/3}\}\rightarrow0$, that is, required
by the Cook's distance. The condition $|c_1|^{-1}n^{2/3}\rightarrow0$
here is to ensure that the influence of the first observation does not
vanish as $n$ goes to infinity. Moreover, in terms of the predictor
dimension $p$, the classical Cook's distance is defined when $p < n$.
Consequently, the condition $n^{-1} p^6 \to0$, or equivalently,
$p=o(n^{1/6})$, constrains the growing rate of $p$ with $n$ at a much
slower rate. We note that although this rate may not be the optimal
one, the condition $p=o(n)$ is clearly necessary for the classical
Cook's distance to be feasible. Next, we examine the condition
$\max\{|c_1|^{-1}(\log n)^{1/2}, l_pp^{-2}c_1^{-4}n\}\rightarrow0$,
that is, required by our new influence measure. For illustration, we
consider a simple case with all the eigenvalues of $\Sigma$ bounded and
$p>n$. We know immediately that both $l_p/p$ and $n/p$ are bounded.
Accordingly, we should have $l_p p^{-2} c_1^{-4} n\rightarrow0$ as long
as $c_1\rightarrow\infty$. As $\log n\rightarrow\infty$ when
$n\rightarrow\infty$, then a sufficient condition for
$\max\{|c_1|^{-1}(\log n)^{1/2},l_p p^{-2} c_1^{-4} n\}\rightarrow0$ is
that $(\log n)^{1/2}/|c_1|\rightarrow0$. This suggests that the
influence point can be consistently detected, as long as $c_1$ diverges
to infinity at a speed faster than $(\log n)^{1/2}$. This is clearly a
rate much slower than $n^{2/3}$. Finally, the bounded eigenvalue
condition \textup{(C.5)} is commonly used in the literature for
estimating covariance matrices \cite{Bick:Levi:2008}. Here it is
assumed for the Cook's distance case. For the new diagnosis measure,
\textup{(C.2)} is required instead, which is weaker
than~\textup{(C.5)}.

\section{Numerical studies}\label{sec3}

We have carried out an intensive simulation study, along with a
microarray data analysis, to examine the empirical performance of our
proposed high-dimensional influence measure. Since the classical Cook's
distance depends on both leverage points and outliers, in our
simulation study, we consider three different scenarios where there
exist outliers only (model~1), leverage points only (model~2), or mixed
leverage points and outliers (model~3). For the scenarios with leverage
points (models 2~and~3), we further consider sub-scenarios where
important covariates contribute to leverage observations, or noisy
covariates contribute to leverage observations. Below we present the
summary of the analysis.

\subsection{Simulation models}\label{sec3.1}

For all simulations, we set the sample size $n = 100$, and the number
of predictors $p = 1000$. We set 10\% of total observations as
influential, so that $\tilde n = 10$. We consider the model
\[
Y_i=\Xbf_i^\top\bolds{\beta}+\varepsilon_i,\qquad i = 1,\ldots, n,
\]
where $\Xbf_i$ is multivariate normal with
$\operatorname{cov}(X_{ij},X_{ij'})=0.5^{|i-j|}$, $\varepsilon_i$
follows the standard normal distribution, and
$\bolds{\beta}=(3,1.5,0,0,2,0,\ldots,0)^\top$. We simulated $n = 100$
i.i.d. observations from this model. Next, we reset the first $\tilde n
= 10$ data observations as coming from another model,
\[
\tilde Y_i=\tilde{\Xbf}{}_i^\top\tilde{\bolds{\beta}}+\varepsilon_i,\qquad i = 1, \ldots,\tilde n,
\]
where perturbations are to be introduced on the regression coefficient,
the covariates and their combination. In particular, we have considered
three perturbation models of generating influential points.
\begin{longlist}[\textit{Model} 1.]
\item[\textit{Model} 1.] The perturbation was introduced on the { response}.
    That is, for $i = 1, \ldots, \tilde n$, $\tilde{\Xbf}_{i}=\Xbf_i$,
    and $\tilde{\bolds{\beta}}=(3,1.5,\kappa,\kappa,2,\kappa,\ldots,
    \kappa)^\top$. In other words, the influential observations { are
    generated according to $\tilde Y_i=\Xbf_i^\top\bolds{\beta}+\kappa
    Z_i+\varepsilon_i$, where $Z_i=\Xbf_i^\top\gamma$ and
    $\gamma=(0,0,1,1,0,1,1,\ldots,1)^\top$. In this case, the responses
    of the influential observations are contaminated by a random
    perturbation $\kappa Z_i$. } Consequently, the corresponding
    responses admit a different pattern, whereas the predictors of
    influential observations follow the same distribution as the rest.

\item[\textit{Model} 2.] The\vspace*{1pt} perturbation was introduced on the predictors and
    keep the response uncontaminated. That is, for $i = 1, \ldots,
    \tilde n$, $\tilde{Y}_i =Y_i$ and $\tilde X_{ij}=X_{ij}+ 30\kappa
    I_{\{j\in S\}}, j=1,\ldots,p$. In other words, a set $S$ of
    predictors admit a different pattern, and its magnitude is
    controlled by the scalar $\kappa$. We examined three choices of
    $S$: $S_1 = \{1,\ldots, 100\}$, and in this case, the influenced
    predictors overlap with those truly relevant ones $\{1,2,5\}$ in
    $\bolds{\beta}$; $S_2 = \{p-100,\ldots, p\}$, and as such there is
    no overlap; and $S_3 = \{1,\ldots, p\}$, and in this case, all
    predictors are subjected to potential influence.

\item[\textit{Model} 3.] The perturbation was introduced on both the
    response and the predictors. That is,
    $\tilde{\bolds{\beta}}=(3,1.5,\kappa,\kappa,2,\kappa,\ldots,
    \kappa)^\top$ and $\tilde X_{ij}=X_{ij}+ 30\kappa I_{\{j\in S\}}$,
    \mbox{$j=1,\ldots,p$}. Again, we considered three sets of $S$ as described
    earlier.
\end{longlist}
It is clear that $\kappa$ is the parameter that dictates the magnitude
of the influential points. When $\kappa=0$, there is no influential
point. We used $\kappa=0$, 0.4, 0.8, 1.2 and~1.6 in our experiment.

\subsection{Performance evaluation}\label{sec3.2}

We evaluate and compare our proposed influence measure in several ways.
First, we study the potential impact of influential data and how the
proposed diagnosis measure could help limit such impact. Toward that
end, we first applied the LASSO or SIS to the full data. Then we
computed the proposed high-dimensional influence measure, evaluated the
corresponding $p$-value, and applied the multiple testing procedure of
\cite{Benjamini:Hochberg:1995}, with the false discovery rate fixed at
$\alpha=$ 5\%. We then obtained a reduced data set by removing those
flagged influential points and applied the LASSO or SIS to the reduced
data set. We evaluated the impact of influential data in terms of
coefficient estimation, variable selection, and variable screening. For
coefficient estimation, we report the error between the estimated and
true $\bolds{\beta}$,
$\mathrm{ERR}=\|\hat{\bolds{\beta}}-\bolds{\beta}_{\mathrm{true}}\|_2$;
for variable selection, we report the false positive rate,
$\mathrm{FPR}=\#\mbox{False Positive}/\#\mbox{True Negative}$; and for
variable screening, we report the coverage probability $\mathrm{CP}$.
In addition, we also report the empirical power of our influence
identification procedure.

Second, we compare our method to two potentially competing solutions in
high-dimensional influence diagnosis. One is a modified Cook's distance
based on the LASSO. That is, we continue to employ the classical Cook's
distance, but estimate the regression coefficient $\bolds{\beta}$ under
a LASSO penalty and as such avoid the difficulty of the OLS estimate
when $p > n$. This seems a very natural solution. We compare it with
our proposal in terms of estimation accuracy, selection accuracy and
power. On the other hand, we note the lack of asymptotic theory for
this modified Cook's distance. To determine the threshold for
influential data, one may use bootstrap. However, in our comparison, we
simply label the observations with the largest $\tilde n$ modified
Cook's distance as influential. This is not feasible in practice, but
provides a useful benchmark for comparison. The other competing
solution is the penalized least absolute deviation via the LASSO
penalty (LAD${}+{}$LASSO)
\cite{Wang:Li:Jiang:2007,Belloni:Chernozhukov:2011}. Due to the use of
the least absolute deviation as the loss function, this method is
designed to handle heavy tailed errors in linear regression, and as
such a potentially useful way to limit impact of the influence
observations.
%
%
\begin{table}
\tabcolsep=0pt
\caption{Simulation results for perturbation model~1. HIM
denotes our proposed high-dimensional diagnosis measure, and CD denotes
the classical Cook's distance}\label{tab:model1}
\begin{tabular*}{\tablewidth}{@{\extracolsep{\fill}}lcd{1.3}d{1.3}d{1.3}d{2.3}d{2.3}@{}}
\hline
&& \multicolumn{5}{c@{}}{$\bolds{\kappa}$}\\[-4pt]
&& \multicolumn{5}{c@{}}{\hrulefill}
\\
\textbf{Method} & \textbf{Criterion} & \multicolumn{1}{c}{\textbf{0}} & \multicolumn{1}{c}{\textbf{0.4}}
& \multicolumn{1}{c}{\textbf{0.8}} & \multicolumn{1}{c}{\textbf{1.2}} & \multicolumn{1}{c@{}}{\textbf{1.6}}
\\
\hline
SIS & CP & 1 &0.25 &0 &0 &0
\\[3pt]
SIS${}+{}$HIM & CP & 1 &1 &1 &1 &1
\\[6pt]
LASSO & ERR & 0.510 & 4.917 & 9.553 & 14.636 &18.478 \\
& FPR & 0.002 & 0.094 &0.103 & 0.107 & 0.106
\\[3pt]
LASSO${}+{}$HIM &ERR& 0.519 & 1.296 & 1.020 & 0.872 & 0.769 \\
&FPR& 0.002 &0.045 &0.029 & 0.015 & 0.012 \\
&Power & \multicolumn{1}{c}{--} & 0.6 & 0.765 & 0.865 & 0.865
\\[6pt]
LASSO${}+{}$CD &ERR& 0.535 & 1.136 & 2.176 & 2.565 & 4.182 \\
&FPR& 0.003 &0.034 &0.066 & 0.072 & 0.076 \\
&Power & \multicolumn{1}{c}{--} & 0.630 & 0.670 & 0.700 & 0.660
\\[3pt]
LAD${}+{}$LASSO & ERR & 0.642 &1.920& 2.073& 2.406& 1.769\\
\hline
\end{tabular*}
\end{table}
%
%
\begin{table}
\tabcolsep=0pt
\caption{Simulation results for perturbation model~\textup{2}. HIM
denotes our proposed high-dimensional diagnosis measure, and CD denotes
the classical Cook's distance}\label{tab:model2}
\begin{tabular*}{\tablewidth}{@{\extracolsep{\fill}}lccd{1.3}d{1.3}d{1.3}d{1.3}d{1.3}@{}}
\hline
&&& \multicolumn{5}{c@{}}{$\bolds{\kappa}$}\\[-4pt]
&&& \multicolumn{5}{c@{}}{\hrulefill}\\
\textbf{Subset} & \textbf{Method} & \textbf{Criterion}& \multicolumn{1}{c}{\textbf{0}} &
\multicolumn{1}{c}{\textbf{0.4}} & \multicolumn{1}{c}{\textbf{0.8}} & \multicolumn{1}{c}{\textbf{1.2}}
& \multicolumn{1}{c@{}}{\textbf{1.6}}
\\
\hline
$S_1$ &SIS & CP & 1 & 0.05 &0 &0 &0 \\
&SIS${}+{}$HIM & CP &1 & 0.05 &0.1 &0.3 &0.25
\\[3pt]
&LASSO &ERR & 0.439 & 4.917 & 4.972 & 4.971 & 4.954 \\
& &FPR & 0.002 & 0.086 & 0.090 & 0.089 & 0.089 \\
& LASSO${}+{}$HIM &ERR& 0.455 & 4.803 & 4.591 & 3.055 & 3.136 \\
& &FPR& 0.002 & 0.080 & 0.060 & 0.055 & 0.044 \\
& &Power & \multicolumn{1}{c}{--} & 0.620 & 0.775 & 0.892 & 0.930
\\[3pt]
& LASSO${}+{}$CD &ERR& 0.513 & 4.566 & 4.568 & 4.603 & 4.533 \\
& &FPR& 0.004 & 0.073 & 0.073 & 0.070 & 0.070 \\
& &Power & \multicolumn{1}{c}{--} & 0.095 & 0.085 & 0.105 & 0.115 \\
& LAD${}+{}$LASSO & ERR & 0.642 & 1.339 & 1.303& 1.320& 1.330
\\[6pt]
$S_2$ &SIS & CP & 1 & 1 &1 &1 &1 \\
&SIS${}+{}$HIM & CP &1 & 1 &1 &1 &1
\\[3pt]
& LASSO &ERR & 0.509 & 0.456 & 0.439 & 0.450 & 0.469 \\
& &FPR & 0.001 &0.001 &0.001 & 0.002 &0.002 \\
& LASSO${}+{}$HIM &ERR & 0.521 & 0.494 & 0.493 & 0.494 & 0.506 \\
& &FPR& 0.001 &0.001 &0.001 & 0.002 &0.002 \\
& &Power & \multicolumn{1}{c}{--} & 0.695 & 0.8 & 0.85 & 0.895
\\[3pt]
& LASSO${}+{}$CD &ERR & 0.548 & 0.523 & 0.532 & 0.556 & 0.551 \\
& &FPR& 0.001 &0.002 &0.002 & 0.002 &0.002 \\
& & Power & \multicolumn{1}{c}{--} & 0.065 & 0.085 & 0.135 & 0.115 \\
&LAD${}+{}$LASSO&ERR& 0.642& 0.650& 0.645& 0.647& 0.634
\\[6pt]
$S_3$ &SIS & CP & 1 &0.35 &0.45 &0.30 &0.25 \\
&SIS${}+{}$HIM & CP &1 &0.50 &0.60 &0.62 &0.65
\\[3pt]
& LASSO &ERR & 0.473 & 1.567 & 1.545 & 1.598 & 1.609 \\
& &FPR & 0.003 &0.051 &0.053 &0.051 & 0.055 \\
& LASSO${}+{}$HIM &ERR & 0.490 & 1.517 & 1.456 & 1.221 & 1.115 \\
& &FPR& 0.003 &0.034 &0.031 &0.023 & 0.033 \\
& &Power & \multicolumn{1}{c}{--} & 0.735 & 0.86 & 0.95 & 0.95
\\[3pt]
& LASSO${}+{}$CD &ERR& 0.560 & 1.751 & 1.700 & 1.743 & 1.871 \\
& &FPR & 0.003 &0.047 & 0.042 & 0.042 & 0.048 \\
& &Power & \multicolumn{1}{c}{--} & 0.115 & 0.085 & 0.115 & 0.110 \\
&LAD${}+{}$LASSO&ERR&0.642& 0.608& 0.573& 0.580&0.581 \\
\hline
\end{tabular*}
\end{table}
%
%
\begin{table}
\tabcolsep=0pt
\caption{Simulation results for perturbation model~\textup{3}. HIM
denotes our proposed high-dimensional diagnosis measure, and CD denotes
the classical Cook's distance}\label{tab:model3}
\begin{tabular*}{\tablewidth}{@{\extracolsep{\fill}}lccd{1.3}d{1.3}d{1.3}d{2.3}d{2.3}@{}}
\hline
&&& \multicolumn{5}{c@{}}{$\bolds{\kappa}$}\\[-4pt]
&&& \multicolumn{5}{c@{}}{\hrulefill}\\
\textbf{Subset} & \textbf{Method} & \textbf{Criterion}& \multicolumn{1}{c}{\textbf{0}}
& \multicolumn{1}{c}{\textbf{0.4}} & \multicolumn{1}{c}{\textbf{0.8}}
& \multicolumn{1}{c}{\textbf{1.2}} & \multicolumn{1}{c@{}}{\textbf{1.6}}
\\
\hline
$S_1$ &SIS & CP & 1 &1 &0.65 &0.10 &0.05 \\
&SIS${}+{}$HIM &CP & 1 &0.90 &1 &1 &1
\\[3pt]
& LASSO &ERR & 0.446 & 1.559 & 5.308 & 9.628 & 14.498 \\
& &FPR & 0.002 &0.062 & 0.093 & 0.099 & 0.098 \\
& LASSO${}+{}$HIM &ERR & 0.447 & 1.278 & 0.771 & 0.499 & 0.542 \\
& &FPR& 0.002 &0.046 & 0.027 & 0.003 & 0.002 \\
& &Power & \multicolumn{1}{c}{--} & 0.185 & 0.94 & 1 & 1
\\[3pt]
& LASSO${}+{}$CD &ERR & 0.559 & 0.686 & 2.149 & 5.623 & 10.926 \\
& &FPR& 0.002 &0.009 & 0.063 & 0.084 & 0.090 \\
& &Power & \multicolumn{1}{c}{--} & 0.555 & 0.720 & 0.675 & 0.585 \\
&LAD${}+{}$LASSO&ERR&0.642& 1.416& 4.367& 8.740& 13.252
\\[6pt]
$S_2$ &SIS & CP & 1 &1 &0.05 &0 &0 \\
&SIS${}+{}$HIM & CP &1 &1 &1 &1 &1
\\[3pt]
&LASSO &ERR & 0.479 & 2.090 & 6.619 & 11.997 & 17.279 \\
& &FPR & 0.002 &0.072 &0.095 & 0.101 &0.101 \\
& LASSO${}+{}$HIM &ERR & 0.494 & 1.836 & 0.696 & 0.475 & 0.501 \\
& &FPR& 0.002 &0.062 &0.009 & 0.002 &0.002 \\
& &Power & \multicolumn{1}{c}{--} & 0.145 & 0.955 & 1 & 1
\\[3pt]
& LASSO${}+{}$CD &ERR & 0.501 & 0.769 & 3.702 & 7.676 & 14.585 \\
& &FPR& 0.003 &0.016 &0.078 & 0.087 &0.091 \\
& & Power & \multicolumn{1}{c}{--} & 0.605 & 0.680 & 0.685 & 0.520 \\
&LAD${}+{}$LASSO&ERR& 0.642& 1.859& 5.855& 10.829& 16.157
\\[6pt]
$S_3$ &SIS & CP & 1 &1 &0.1 &0 &0 \\
&SIS${}+{}$HIM & CP & 1&1 &1 &1 &1
\\[3pt]
& LASSO &ERR & 0.464 & 1.682 & 5.720 & 10.943 & 17.384 \\
& &FPR & 0.002 &0.065 &0.098 & 0.103 & 0.105 \\
& LASSO${}+{}$HIM &ERR & 0.484 & 1.479 & 1.262 & 0.557 & 0.515 \\
& &FPR& 0.002 &0.057 &0.034 & 0.003 & 0.002 \\
& &Power & \multicolumn{1}{c}{--} & 0.1 & 0.87 & 1 & 1
\\[3pt]
& LASSO${}+{}$CD &ERR & 0.586 & 0.726 & 1.874 & 4.504 & 7.566 \\
& &FPR& 0.002 &0.013 &0.055 & 0.074 & 0.087 \\
& & Power & \multicolumn{1}{c}{--} & 0.465 & 0.765 & 0.810 & 0.855 \\
& LAD${}+{}$LASSO &ERR& 0.642& 1.635& 5.264& 10.662& 17.023 \\
\hline
\end{tabular*}
\end{table}

\subsection{The results}\label{sec3.3}
The averages of a total of 200 random replications are reported in
Tables~\ref{tab:model1}--\ref{tab:model3}. We make the following
observations.
\begin{longlist}[(3)]
\item[(1)] First, the presence of influential points significantly
    affects variable selection and screening accuracy.\vadjust{\goodbreak} This can be seen
    by comparing the results between SIS and SIS${}+{}$HIM in terms of CP.
    Consider, for example, Table~\ref{tab:model1}. As $\kappa$ increases, the coverage
    probability of the SIS method deteriorates quickly from 1 with
    $\kappa=0$ to 0 with $\kappa=1.6$. This confirms that influential
    observations do affect variable screening consistency. Meanwhile,
    the performance of SIS${}+{}$HIM is quite encouraging as its CP values
    maintains at 1 for every $\kappa$ value considered. This suggests
    that the proposed HIM method helps SIS in removing the influential
    observations.

\item[(2)] Second, the presence of influential observations does
    affect estimation accuracy seriously. This can be seen clearly
    by comparing the results of LASSO and LASSO${}+{}$HIM in terms of ERR
    values. For instance, the ERR values in Table~\ref{tab:model3} for LASSO with
    $S_1$ increases quickly from 0.446 with $\kappa=0$ to 14.498
    with $\kappa=1.6$. This confirms that influential observations
    do affect the accuracy of the LASSO estimate in a negative way.
    However, we find that the ERR values of LASSO${}+{}$HIM are always
    well controlled with ERR${}<{}$1.5. In fact, as $\kappa$
    increases, the power for HIM to detect influential observation
    increases. Thus, those influential observations are more likely
    to be detected and eliminated from the data analysis. This
    makes the ERR values of LASSO${}+{}$HIM eventually converges to a
    level around ERR${}\approx{}$0.5, as $\kappa$ increases. This
    confirms the usefulness of the HIM method for LASSO estimation,
    even though its definition only involves marginal correlation
    coefficients.

\item[(3)] Third, the performance of LASSO${}+{}$CD is mixed. If the
    perturbation is due to the response only as in Table~\ref{tab:model1}, it does
    yield much better performance than LASSO with much smaller ERR
    values. This suggests that LASSO${}+{}$CD can perform well to limit the
    effect of influential points. However, even for this example, it is
    still outperformed by LASSO${}+{}$HIM. However, the story changes if the
    perturbation is due to the predictors as in Table~\ref{tab:model2}. This is to be
    expected because, with contaminated predictors, LASSO is no longer
    a stable method for variable selection. If predictors are selected
    incorrectly, the subsequent modified Cook's distance cannot be
    calculated appropriately. This makes the performance of LASSO${}+{}$CD
    unsatisfactory.

\item[(4)] Fourth, as a robust regression method, we find that
    LAD${}+{}$LASSO performs quite well. Its ERR values are smaller than
    those of the LASSO estimates in all the tables. However, in most
    cases, it is still outperformed by LASSO${}+{}$HIM as seen
    from Tables~\ref{tab:model1}~and~\ref{tab:model3}.

\item[(5)] Lastly, we find that for most cases, the reported FPR
    values are well controlled. Furthermore, as $\kappa$ increases,
    the corresponding empirical power increases toward 100\%. These
    findings are consistent with the theoretical claims in Theorems
    \ref{Dk}~and~\ref{Cook-u}.
\end{longlist}
To summarize, our simulation experiments confirm that the proposed HIM
method is useful in controlling the effects of the influential
observations in terms of parameter estimation and variable screening.

\subsection{A real data example}\label{sec3.4}

We applied our proposed influence diagnosis approach to a microarray
data of \cite{Scheetz:etal:2006}, and noted that the analysis results
become substantially different when the detected influential
observations are removed. For this dataset, F1 animals were
intercrossed and then 120 twelve-week-old male offspring were selected
for tissue harvesting from the eyes and for microarray analysis.
The Affymetrix microarrays that were used to analyze the RNA from the eyes
of those F2 animals contain over 31,042 different probe sets. Among
them, one probe is for gene TRIM32, which was recently found to cause
Bardet--Biedl syndrome \cite{Chiang:etal:2006}, a genetically
heterogeneous disease of multiple organ systems including the retina.
One goal of interest of this data analysis is to find genes whose
expressions are correlated with that of gene TRIM32. We first followed
\cite{Huang:Ma:Zhang:2008} to exclude probes that were not expressed in
the eye or that lacked sufficient variation, which results in 18,975
probes as regressors. We then followed \cite{Fan:Lv:2008} to retain the
top 1000 probes that are mostly correlated with the probe of TRIM32.
The resulting analysis has $p = 1000$ predictors and a sample size $n =
120$. As a standard procedure \cite{Huang:Ma:Zhang:2008}, all the
probes are standardized to have mean zero and standard deviation one.
%
%
\begin{figure}

\includegraphics{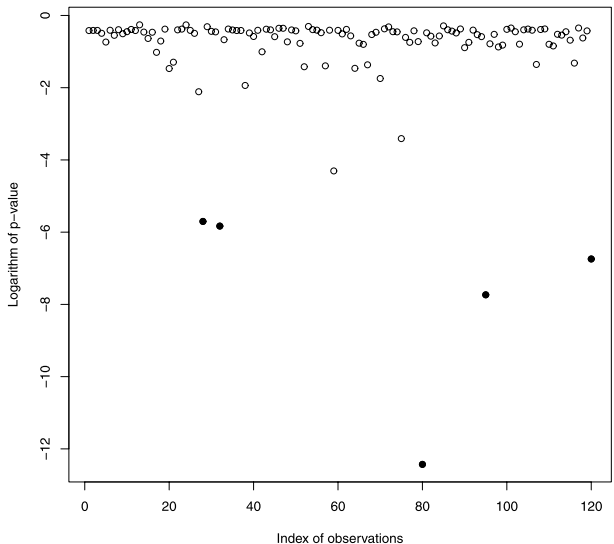}

\caption{The logarithm of the $p$-value for each observation:
the detected influential points are denoted by solid circles.}\label{pvalue1}
\end{figure}

We next applied our method with FDR rate $\alpha=0.10$ to the data, and
identified a total of 5 influential observations. Their corresponding
$p$-values were 0, 0.0004, 0.0011, 0.0029 and 0.0033, respectively. We
also show the logarithm of $p$-values versus the indices for all 120
observations in Figure~\ref{pvalue1}. To assess the influence of the
detected points, we again compared the LASSO estimate with and without
those points. 
Since we used ten-fold cross-validation to select the tuning parameter
and every run is random, we repeated this analysis 100 times and report
the average results.

We summarize the difference of the estimates in the following aspects:
the sparsity, the norm difference and the angle between the two
estimates. First, by removing the identified influential observations,
the resulting LASSO estimate is considerably more sparse. The average
model size with the full data is~63. By contrast, the average model
size without the influential observations reduces to 27 on the average.
The existence of the potential influential points clearly shows a
significant effect on the model size. Besides that, the average number
of the commonly selected predictors by fitting the full data and the
reduced data, respectively, is only 8.67, which again shows clear
discrepancy of the two estimates. Consequently, the influential points
identified by our approach seem to have significant effect for
subsequent analysis.
Second,\vspace*{-1pt} denote $d_{0}=\|\hat{\bolds{\beta}}_{\mathrm{full}}\|_2$,
$d_{1}=\|\hat{\bolds{\beta}}_{\mathrm{redu}}\|_2$ and
$d_{2}=\|\hat{\bolds{\beta}}_{\mathrm{redu}}-\hat{\bolds{\beta}}_{\mathrm{full}}\|_2$,
where $\hat{\bolds{\beta}}_{\mathrm{full}}$ is the\vspace*{1pt} LASSO estimate using
all the observations and $\hat{\bolds{\beta}}_{\mathrm{redu}}$ is the
estimate after removing the influential points identified by~HIM. We
observe that the average of $(d_{0}-d_{1})/d_{0}$ is 0.532 and that of
$d_{2}/d_{0}$ is 0.972. Both show that the estimates without
influential points are quite different in terms of the $\ell_2$ norm.
In addition, the angle\vspace*{-1pt} between $\hat{\bolds{\beta}}_{\mathrm{full}}$
and $\hat{\bolds{\beta} }_{\mathrm{redu}}$, which is defined as
$\hat{\bolds{\beta}}{}_{\mathrm{full}}^\top\hat{\bolds{\beta}
}_{\mathrm{redu}}/d_{0}d_{1}$, equals 0.262, averaged over $100$ times.
These numbers again indicate that the estimates change substantially
after removing the influential observations. In summary, this analysis
illustrates the importance of influence diagnosis, and the identified
influential observations should be treated with extreme care.

\section{Extension to generalized linear models}\label{sec4}

The main idea of the high-dimensional influence measure can be extended
to a broad class of regression models. Here we briefly discuss one such
extension to generalized linear models (GLM). Assume that the data
$(\Xbf_i, Y_i)$ follow an exponential distribution with the canonical
probability density function, $f(y;\theta)=\exp\{y\theta- b(\theta) +
c(y)\}$, and the conditional mean is of the form
\[
E(Y_i\mid \Xbf_i)=b'\bigl(\theta(\Xbf_i)\bigr)=g^{-1} \bigl( \beta_0 +\Xbf_i^\top\bolds{\beta}_1 \bigr),
\]
where $g$ is a known link function. For the purpose of feature
screening in ultra high-dimensional regressions, \cite{Fan:Song:2010}
introduced a marginal utility measure, the maximum marginal likelihood
estimator, as
\[
\hat{\bolds{\beta}}_j=(\hat{\beta}_{j,0}, \hat{\beta}_j)=\arg\min E_n l(Y,\beta_{j0}+\beta_j X_j),
\]
where $l(Y;\theta)=-Y\theta+b(\theta)+\log c(Y)]$ and $E_n
f(X,Y)=n^{-1}\sum_{i=1}^n f(X_i,Y_i)$. That is, $\hat{\bolds{\beta}}_j$
is the maximum likelihood estimator of fitting a GLM model of $Y$ on
the $j$th predictor $X_j$ alone plus an intercept. As remarked by
\cite{Fan:Song:2010}, this measure can be rapidly computed.

In the context of high-dimensional diagnosis, we define the
high-dimensional influence measure for generalized linear models for
the $k$th observation, \mbox{$k = 1, \ldots, n$,} as
%
%
%
\begin{eqnarray}\label{hcd-glm}
\mathcal{D}_k^{\mathrm{glm}} &=& \frac{1}{p}\sum_{j=1}^p \bigl\llVert\hat{\bolds{\beta}}_j-\hat{\bolds{\beta}}{}_{j}^{(k)} \bigr\rrVert_2^2,
\end{eqnarray}
where $\hat{\bolds{\beta}}{}_j^{(k)}$ denotes\vspace*{-3pt} the maximum marginal
likelihood estimator but with the $k$th observation removed. For GLM,
the estimator $\hat{\bolds{\beta}}_j$ and
$\hat{\bolds{\beta}}{}_{j}^{(k)}$ may not\vspace*{-2pt} have a closed-form solution.
Consequently, the exact distribution of the proposed statistic
$\mathcal{D}_k^{\mathrm{glm}}$ is complicated and some approximation is
necessary. The detailed derivation, however, is beyond the scope of
this paper. In practice, one can always sort the values of
$\{\mathcal{D}_k^{\mathrm{glm}}, k=1,\ldots,n\}$ and remove those
observations associated with large values~of~$\mathcal{D}_k^{\mathrm{glm}}$.

We have conducted a small simulation study to examine the empirical
performance of this measure for GLM. The simulation setup is similar to
that of model~1 in Section~\ref{sec3.1}, except that this time we adopt
a binary response model, $P(Y_i=1\mid \Xbf_i)=1/[1+\exp\{-(2 +
\Xbf_i^\top \bolds{\beta})\}]$, where
$\bolds{\beta}=\bolds{\beta}_{\mathrm{true}}=(5,5,0,\ldots,0)^\top$,
and the outliers are generated from the model
$\bolds{\beta}=\bolds{\beta}_{\mathrm{infl}}=(5,5,0,\ldots,0,-\kappa,
\ldots, -\kappa)^\top$ with $p/2$ many $\kappa$'s. We set $n=100$, with
10\% influential observations, that is, $n_{\mathrm{infl}}=10$, and we set
$p = 50$ or 100. Since the asymptotic distribution of
$\mathcal{D}_k^{\mathrm{glm}}$ is not available for the logistic
regression, we flag the 10 observations with the largest $p$-values of
$\mathcal{D}_k^{\mathrm{glm}}$ as influential. For a binary response,
one is often interested in classification. As such we compare the
misclassification error rate for the full data as $E_{\mathrm{full}}$
and for the reduced data as $E_{\mathrm{redu}}$ without the detected
influential points. We also report the empirical power. The results out
of 200 data replications are summarized in Table~\ref{powsiz_md3}.
%
\begin{table}
\tabcolsep=0pt
\tablewidth=250pt
\caption{Simulation results for the logistic model}\label{powsiz_md3}
\begin{tabular*}{250pt}{@{\extracolsep{\fill}}lcccccc@{}}
\hline
&& \multicolumn{5}{c@{}}{$\bolds{\kappa}$} \\[-4pt]
&& \multicolumn{5}{c@{}}{\hrulefill}
\\
$\bolds{p}$ & \textbf{Criterion} & \textbf{0} & \textbf{0.4} & \textbf{0.8} & \textbf{1.2} & \textbf{1.6}\\
\hline
\phantom{0}50 &Power & \multicolumn{1}{c}{--} & 0.220 & 0.472 & 0.422 & 0.254 \\
&$E_{\mathrm{full}}$ & 0.037 & 0.062 & 0.088 & 0.083 & 0.064 \\
&$E_{\mathrm{redu}}$ & 0.018 & 0.022 & 0.031 & 0.049 & 0.033
\\[6pt]
100 &Power & \multicolumn{1}{c}{--} & 0.332 & 0.386 & 0.220 & 0.152 \\
&$E_{\mathrm{full}}$ &0.047 & 0.069 & 0.065 & 0.045 & 0.029 \\
& $E_{\mathrm{redu}}$ &0.020 & 0.042 & 0.028 & 0.018 & 0.019 \\
\hline
\end{tabular*}
\end{table}

From Table~\ref{powsiz_md3}, we note that the proposed method has some
power for a logistic model, but it is lower than that in a linear
model. This is probably due to the fact that a binary response contains
much less information, and thus detecting influential observations in a
logistic model is much more challenging, especially in a
high-dimensional setting. On the other hand, we also observe from
Table~\ref{powsiz_md3} that removing those points with the largest
values of $\mathcal{D}_k^{\mathrm{glm}}$ improves the classification
accuracy by a large margin. This again suggests that the usefulness of
influence diagnosis. Meanwhile, we remark that further investigation
into both theoretical and empirical properties of high-dimensional
influence measure in GLM is warranted.

\section{Conclusion}\label{sec5}

We perceive several future avenues to extend the proposed work in this
article. First, we have employed the leave-one-out principle when
quantifying influence of individual observations. We expect that our
high-dimensional influence measure can also be generalized to the cases
of leaving out pairs of observations, or triplets or more. Such a
strategy can be useful when those observations conceal one another
\cite{Draper:Smith:1998}. Second, we have focused on the classical
linear model in our development, while extension to more sophisticated
models, such as the generalized linear model, that is, briefly examined
in Section~\ref{sec4}, survival models, and semiparametric additive
models, deserve further investigations. Finally, our proposal deals
with the cross sectional data with i.i.d. observations. It is
interesting to extend the proposed influence measure to complex
correlated data such as longitudinal data where dependence among
observations needs to be taken into consideration in influence
diagnosis \cite{Zhu:Ibrahim:Cho:2012}.

\begin{appendix}\label{sec6}
\section*{Appendix}

We outline the main idea of the proof for the asymptotic distribution
of $\mathcal{D}_k$ in Theorem~\ref{Dk}. First, we decompose
$\mathcal{D}_k$ as $\mathcal{D}_k = B_1+B_2+B_3-2B_4$ as given in
Section~\ref{sec2.3}.\vadjust{\goodbreak} Then we compute the mean and variance of $B_i$,
$i=1,\ldots,4$ as presented in Proposition~\ref{i-cook}. This step
builds on the assumption of normality of the predictors and benefits
from the fact that the predictor dimension goes to infinity. Comparing
the orders of the variance of $B_i$, we find that $B_2$ is the leading
term. We then study the asymptotic distribution of $B_2$,
which turns out to follow a $\chi^2(1)$ distribution. Recall in Section~\ref{sec2.3}, we defined $K_{p,tl}=p^{-1}\sum_{j=1}^p
X_{tj}X_{lj}$, for $t,l=1,\ldots,n$, $l_p = \sum_{j=1}^{p} \lambda_j^2
= O(p^r)$ and $c_p = \max_{1 \leq j \leq p} \lambda_j$, where
$\lambda_j$'s are the eigenvalues of the covariance matrix
$\bolds{\Sigma}$. Furthermore, we define $a_p = \sum_{j=1}^{p}
\lambda_j^4$, $C_1=E(Y_tY_l K_{p,tl})^2$ and $C_2=E[Y_t^2(\sum_{j=1}^p
\rho_jX_{tj}/p)^2]$ for any $t \neq l$.


\begin{pf*}{Proof of Proposition~\ref{i-cook}}
We break the proof into three parts: first, we obtain an expansion of
$\mathcal{D}_k$; second, we derive $E(\mathcal{D}_k)$; and finally, we
derive the asymptotic behavior of the components in the expansion of
$\mathcal{D}_k$.
\begin{longlist}[\textit{Step} 3.]
\item[\textit{Step} 1.] First, we have the following expansion
    for $\mathcal{D}_k$, $k = 1, \ldots, n$:
\begin{eqnarray*}
\mathcal{D}_k&=&\frac{1}{p}\sum_{j=1}^p\Biggl(\frac{1}{n-1}\sum_{t=1,t\neq k}^n Y_tX_{tj}-\frac{1}{n}\sum_{t=1}^n Y_tX_{tj}\Biggr)^2
\\[-0.5pt]
&=& \frac{1}{p}\sum_{j=1}^p\Biggl\{\frac{1}{n(n-1)}\sum_{t=1,t\neq k}^nY_tX_{tj}-\frac{1}{n} Y_kX_{kj}\Biggr\}^2
\\[-0.5pt]
&=&\frac{1}{p}\sum_{j=1}^p \Biggl\{\frac{1}{n(n-1)}\sum_{t=1,t\neq k}^n Y_tX_{tj} \Biggr\}^2 +\frac{1}{pn^2} Y_k^2\sum_{j=1}^p
X_{kj}^2
\\[-0.5pt]
&&{} -\frac{2}{pn^2(n-1)}\sum_{t=1,t\neq k}^n Y_tY_k \Biggl\{\sum_{j=1}^pX_{tj}X_{kj}\Biggr\}
\\[-0.5pt]
&=& \frac{1}{p\{n(n-1)\}^2}\sum_{t=1,t\neq k}^nY_t^2 \Biggl\{ \sum_{j=1}^p
X_{tj}^2 \Biggr\}+\frac{1}{pn^2}Y_k^2\Biggl\{ \sum_{j=1}^p X_{kj}^2\Biggr\}
\\[-0.5pt]
&&{} +\frac{1}{p\{n(n-1)\}^2}\sum_{t\neq s, t,s\neq k}Y_tY_s \Biggl\{ \sum_{j=1}^p
X_{tj} X_{sj} \Biggr\}
\\[-0.5pt]
&&{} -\frac{2}{pn^2(n-1)} \sum_{t=1,t\neq k}^n Y_kY_t \Biggl\{\sum_{j=1}^p X_{tj}X_{kj} \Biggr\}
\\
&=&\frac{1}{\{n(n-1)\}^2}\sum_{t=1}^n Y_t^2 K_{p,tt}+ \biggl[\frac{1}{n^2}-\frac{1}{\{n(n-1)\}^2} \biggr] Y_k^2K_{p,kk}
\\
&&{} +\frac{1}{\{n(n-1)\}^2}\sum_{t\neq s} Y_tY_s K_{p,t s}
\\
&&{}- \biggl[\frac{2}{\{n(n-1)\}^2}+\frac{2}{n^2(n-1)} \biggr] \sum_{t=1,t\neq k}^n
Y_kY_t K_{p,tk}
\\
&=&\frac{1}{\{n(n-1)\}^2}\sum_{t=1}^n Y_t^2 K_{p,tt}+\frac{(n-2)}{n(n-1)^2}Y_k^2K_{p,kk}
\\
&&{} +\frac{1}{\{n(n-1)\}^2} \sum_{t\neq s}Y_tY_s K_{p,t s} -\frac{2}{n(n-1)^2} \sum_{t=1,t\neq k}^n Y_kY_t K_{p,tk}
\\
&:=& B_1+B_2+B_3-2B_4.
\end{eqnarray*}

\item[\textit{Step} 2.] Next, we derive the expectation of $\mathcal{D}_k$.
    It is easy to see that
\begin{eqnarray*}
E(B_1)&=&\frac{1}{pn(n-1)^2} \sum_{j=1}^p E(Y_k^2 X_{kj}^2),
\\
E(B_2)&=&\frac{n-2}{pn(n-1)^2} \sum_{j=1}^p E(Y_k^2 X_{kj}^2),
\\
E(B_3)&=&\frac{1}{pn(n-1)}\sum_{j=1}^p \rho_j^2,\qquad
E(B_4)=\frac{1}{pn(n-1)}\sum_{j=1}^p \rho_j^2.
\end{eqnarray*}
Therefore, we have
\[
E(\mathcal{D}_k)=E(B_1+B_2+B_3-2B_4)=\frac{1}{pn(n-1)} \sum_{j=1}^p\operatorname{var} (Y_k X_{kj}).
\]
By Lemmas 1 and 3, we have
\begin{eqnarray*}
E\bigl\{Y_k^2\bigl(K_{p,kk}-E(K_{p,kk})\bigr)\bigr\}&\leq& E^{1/2}\bigl(Y_k^4\bigr)E^{1/2}
\bigl[\bigl\{ K_{p,kk}-E(K_{p,kk})\bigr\}^2\bigr]
\\
&=&O \bigl(p^{-1} l_p^{1/2}\bigr)
\end{eqnarray*}
and
\[
p^{-1}\sum_{j=1}^pE^2(Y_kX_{kj})=p^{-1}\sum_{j=1}^p \rho_j^2=O\bigl(p^{-1}c_p \bigr).
\]
In addition, noting that $c_p^2\leq l_p $, we have
\begin{eqnarray*}
&& p^{-1}\sum_{j=1}^p \operatorname{var}(Y_kX_{kj})
\\
&&\qquad =p^{-1}\sum_{j=1}^p \bigl\{ E\bigl(Y_k^2X_{kj}^2\bigr)-E^2(Y_kX_{kj}) \bigr\}
\\
&&\qquad = E\bigl(Y_k^2\bigr)E(K_{p,kk}) +E\bigl\{Y_k^2\bigl(K_{p,kk}-E(K_{p,kk})\bigr) \bigr\}-p^{-1}\sum_{j=1}^p\rho^2_j
\\
&&\qquad = E\bigl(Y_k^2\bigr)E(K_{p,kk})+O\bigl(p^{-1} l_p^{1/2}\bigr).
\end{eqnarray*}
Consequently, we have
\[
E(\mathcal{D}_k)=\frac{1}{\{p n(n-1)\}}E\bigl(Y_k^2
\bigr)\sum_{j=1}^p E\bigl(X_{kj}^2\bigr)+O\bigl(n^{-2}p^{-1}l_p^{1/2}
\bigr).
\]

\item[\textit{Step} 3.] Next, we derive the asymptotic behavior of $B_i,
    i=1,\ldots, 4$.

\begin{longlist}[Step 3.1.]
\item[\textit{Step} 3.1.] We start with the variance of  $B_1$.  Note that
\[
\operatorname{var}(B_1)=\frac{n}{n^4(n-1)^4} \operatorname{var}
\bigl(Y_t^2 K_{p,tt}\bigr)
\]
and that $E(Y_t^4 K_{p,tt}^2)\leq E^{1/2}(Y_t^8)E^{1/2}(K_{p,tt}^4)$.
Furthermore,
\begin{eqnarray*}
E(K_{p,tt})^{4} & = & p^{-4} E\bigl(
\Zbf_t^\top\bolds{\Sigma}\Zbf_t
\bigr)^4 =p^{-4} E\Biggl[\sum_{j=1}^p
\lambda_j \bigl(\Zbf_t^\top
\ubf_j\bigr)^2\Biggr]^4
\\
& \leq& p^{-4} E\bigl[\bigl(\Zbf_t^\top
\ubf_j\bigr)^8\bigr] \Biggl(\sum
_{j=1}^p \lambda_j\Biggr)^4
\leq E\bigl(\Zbf_t^\top\ubf_j
\bigr)^8.
\end{eqnarray*}
The last equation holds because $\operatorname{tr}(\Sigma)=\sum_{j=1}^p
\lambda_j=p$. As a result, we have
\[
\operatorname{var}(B_1)=O\bigl(n^{-7}\bigr).
\]

\item[\textit{Step} 3.2.] Next, we consider the variance of  $B_4$. By Lemma~3,
    we have
\[
E(B_4)=\frac{1}{pn(n-1)} \sum_{j=1}^p\rho_j^2=O\bigl(c_p\cdot p^{-1}n^{-2}\bigr).
\]
In addition, it is easy to see
\begin{eqnarray*}
E\bigl(B_4^2\bigr)&=&\frac{1}{n^2(n-1)^4} \Biggl\{\sum
_{t=1, t\neq k}^n E\bigl(Y_k^2Y_t^2 K_{p,tk}^2\bigr)+\sum_{t\neq s}^{ t,s\neq k}
E\bigl(Y_k^2Y_tY_s K_{p,tk}K_{p,sk}\bigr) \Biggr\}
\\
&=&\frac{1}{n^2(n-1)^4} \Biggl[ (n-1) E\bigl(Y_k^2Y_t^2
K_{p,tk}^2\bigr)
\\
&&\hspace*{55pt}{} +(n-1) (n-2)E \Biggl\{Y_k^2
\Biggl(\sum_{j=1}^p \rho
_jX_{kj}/p \Biggr)^2 \Biggr\} \Biggr]
\\
&=&\frac{1}{n^2(n-1)^3}E\bigl(Y_k^2Y_t^2
K_{p,tk}^2\bigr)+ \frac
{(n-2)}{n^2(n-1)^3} E \Biggl
\{Y_k^2 \Biggl(\sum_{j=1}^p
\rho_jX_{kj}/p \Biggr)^2 \Biggr\}
\\
&=&\frac{1}{n^2(n-1)^3} C_1+ \frac{(n-2)}{n^2(n-1)^3}C_2,
\end{eqnarray*}
where $C_1$, $C_2$ are defined as in (1.1) and (1.2) in
the supplementary material \cite{suppJZ}, respectively. From the proof of Lemma 2, we know
$C_2=O(c_p^2 p^{-2})$ and $C_1=C_{11}+C_{12}$, with $C_{11}=O(p^{-2}
a_p^{1/2})$ and $C_{12}=l_p p^{-2} E^2(Y_t^2) $. Therefore, we have
\begin{eqnarray}
E\bigl(B_4^2\bigr)&=&\frac{l_p}{p^2n^2(n-1)^3} E^2
\bigl(Y_t^2\bigr) +O\bigl(n^{-5}p^{-2}a_p^{1/2}
\bigr)+O\bigl(c_p^2 p^{-2}n^{-4}\bigr)
\nonumber
\\
&=& O\bigl( l_p p^{-2}n^{-5}\bigr)+O
\bigl(c_p^2 p^{-2}n^{-4}\bigr).
\nonumber
\end{eqnarray}
Consequently, we have
\[
\operatorname{var}(B_4)= E\bigl(B_4^2
\bigr)-E^2(B_4)=O\bigl(l_p
p^{-2}n^{-5}\bigr)+O\bigl(c_p^2p^{-2}n^{-4}
\bigr).
\]

\item[\textit{Step} 3.3.] Next, we aim at
    $\operatorname{var}(B_3)$. First, we have
\begin{eqnarray*}
B_3&= &\frac{1}{p\{n(n-1)\}^2}\sum_{t\neq s}Y_tY_s
\Biggl(\sum_{j=1}^pX_{sj}X_{tj}\Biggr)
\\
&=&\frac{1}{pn(n-1)} \sum_{t\neq s}\phi(Y_t,Y_s, \Xbf_s, \Xbf_t)/\bigl\{(n-1)n\bigr\}:=\frac{1}{pn(n-1)} \bar{B}_3,
\end{eqnarray*}
where $ \phi(Y_t,Y_s,\Xbf_s,\Xbf_t)=\sum_{t\neq
s}Y_tY_s(\sum_{j=1}^p X_{sj}X_{tj})$. Let
\begin{eqnarray*}
\phi_1(Y_t, \Xbf_t )&=&E \bigl
\{\phi(Y_t,Y_s,\Xbf_s,\Xbf_t)-
E(\phi(Y_t,Y_s,\Xbf_s,
\Xbf_t)\mid Y_t, \Xbf_t \bigr\}
\\
&=&Y_t\sum_{j=1}^p
X_{tj}\rho_j -\sum_{j=1}^p
\rho_j^2.
\end{eqnarray*}
Noting that $\bar{B}_3$ is an $U$-statistic, and by the properties of
the $U$-statistic, we have
%
\begin{eqnarray}\label{var(A2)}
\operatorname{var}(B_3)&=&\frac{1}{\{pn(n-1)\}^2}\operatorname{var}(\bar{B}_3)\nonumber
\\
&=&\frac{1}{\{pn(n-1)\}^2} \biggl[\frac{4}{n} \operatorname{var}\bigl\{\phi_1(Y_t,
\Xbf_t)\bigr\}+o\bigl(n^{-2}\bigr) \biggr]
\\
&=&\frac{4}{p^{2}n^3(n-1)^2} \operatorname{var} \Biggl\{Y_t \sum
_{j=1}^p X_{tj}\rho_j \Biggr\}+o\bigl(p^{-2}n^{-4}(n-1)^{-2}\bigr).\nonumber
\end{eqnarray}
Furthermore, we have
%
%
%
\begin{eqnarray}\label{var_YXR}
&& p^{-2} \operatorname{var}
\Biggl(Y_t \sum_{j=1}^p
X_{tj}\rho_j \Biggr)\nonumber
\\[-2pt]
&&\qquad =E \Biggl\{Y_t^2 \Biggl(\sum_{j=1}^p X_{tj} \rho_j/p \Biggr)^2 \Biggr\} - \Biggl(\sum_{j=1}^p \rho_j^2/p\Biggr)^2
\\[-2pt]
&&\qquad \leq E^{1/2}\bigl(Y_t^4\bigr)E^{1/2}
\Biggl\{ \Biggl(\sum_{j=1}^p X_{tj}\rho_j/p \Biggr)^4 \Biggr\} - \Biggl(
\sum_{j=1}^p \rho_j^2/p \Biggr)^2.\nonumber
\end{eqnarray}
In addition, we show in the proof of Lemma 2 that
%
%
\begin{equation}\label{Xrho4}
E \Biggl\{ \Biggl(\sum_{j=1}^p X_{tj}\rho_j/p \Biggr)^4 \Biggr\} =O\bigl(c_p^4 p^{-4}\bigr)
\end{equation}
and in Lemma 3 that
%
\begin{equation}\label{rho2}
\sum_{j=1}^p
\rho_j^2/p=O\bigl(c_p\cdot p^{-1}
\bigr).
\end{equation}
Combining (\ref{var(A2)})--(\ref{rho2}), we have
\[
\operatorname{var}(B_3)=O\bigl(c_p^2
n^{-5} p^{-2}\bigr)+O\bigl(p^{-2}n^{-6}
\bigr).
\]

\item[\textit{Step} 3.4.] Finally we turn to $B_2$, which can
    be written as
\begin{eqnarray*}
B_2&=&\frac{(n-2)}{n(n-1)^2}Y_k^2K_{p,kk}
\\[-2pt]
&=& \frac{(n-2)}{n(n-1)^2}
\bigl[Y_k^2 \bigl\{K_{p,kk}-E(K_{p,kk}) \bigr\}+Y_k^2 E(K_{p,kk})
\bigr]
\\[-2pt]
&:=&B_{21}+B_{22}.
\end{eqnarray*}
By Lemma 1, we have
\begin{eqnarray*}
E\bigl(Y_k^4\bigl\{K_{p,kk}-E(K_{p,kk})
\bigr\}^2\bigr) & \leq& E^{1/2}\bigl(Y_k^8
\bigr)E^{1/2}\bigl(\bigl\{K_{p,kk}-E(K_{p,kk})\bigr\}^4\bigr)
\\[-2pt]
&=&O\bigl(p^{-2} l_p\bigr),
\\[-2pt]
E\bigl(Y_k^2\bigl\{K_{p,kk}-E(K_{p,kk})
\bigr\}\bigr) & \leq& E^{1/2}\bigl(Y_k^4
\bigr)E^{1/2}\bigl(\bigl\{K_{p,kk}-E(K_{p,kk})\bigr\}^2\bigr)
\\[-2pt]
&=&O\bigl(l_p^{1/2} p^{-1}\bigr).
\end{eqnarray*}
Therefore, we have
\begin{eqnarray*}
\operatorname{var}(B_{21}) &=& \biggl\{\frac{(n-2)}{n(n-1)^2} \biggr
\}^2 \operatorname{var}\bigl(Y_k^2
\bigl[K_{p,kk}-E(K_{p,kk})\bigr]\bigr)=O\bigl(n^{-4}p^{-2}
l_p\bigr),
\\[-2pt]
\operatorname{var}(B_{22}) &=& O\bigl(n^{-4}\bigr).
\end{eqnarray*}
This completes the proof.\quad\qed
\end{longlist}
\end{longlist}\noqed
\end{pf*}\eject

\begin{pf*}{Proof of Theorem~\ref{Dk}}
Consider the behavior of $K_{p,kk}$, $k = 1, \ldots, n$, for a~sufficient large $p$
\[
K_{p,kk}=\sum_{j=1}^p X_{kj}^2/p=\Xbf_k^\top
\Xbf_k/p=\Zbf_k^\top\bolds{\Sigma} \Zbf_k=\sum_{j=1}^p \lambda_j \bigl(\Zbf_k^\top\ubf_j \bigr)^2/p.
\]
Its variance is $ \operatorname{var}(K_{p,kk})=2\sum_{j=1}^p
\lambda_j^2/p^2=2p^{-2} l_p$. Under the\vspace*{-1pt} assumption $ l_p=O(p^r)$ with
$0\leq r<2$, we have $K_{p,kk}= E(K_{p,kk}) +O_p(p^{r/2-1})$, and
consequently,
\[
Y_k^2K_{p,kk}=Y_k^2\bigl[ E(K_{p,kk})+O_p\bigl(p^{r/2-1}\bigr)\bigr].
\]
In addition, noting that $E[Y_k^2(K_{p,kk}- E(K_{p,kk}))]\leq
E^{1/2}(Y_k^4)(\operatorname{var}(K_{p,kk}))^{1/2}=O(p^{r/2-1})$, we
have
\begin{eqnarray*}
E\bigl(Y_k^2K_{p,kk}\bigr)&=&E
\bigl(Y_k^2\bigr) E(K_{p,kk})+E
\bigl[Y_k^2\bigl(K_{p,kk}- E(K_{p,kk})
\bigr)\bigr]
\\
&=&E\bigl(Y_k^2\bigr) E(K_{p,kk})+O
\bigl(p^{r/2-1}\bigr).
\end{eqnarray*}
Therefore, we have
\[
Y_k^2K_{p,kk}-E\bigl(Y_k^2K_{p,kk}
\bigr)=\bigl[Y_k^2-E\bigl(Y_k^2
\bigr)\bigr] E(K_{p,kk})+O_p\bigl(p^{r/2-1}\bigr).
\]
As a result, it holds that%
%
\begin{equation}
\label{B2} B_2-E(B_2)=\frac{n-2}{n(n-1)^2} \bigl\{
\bigl[Y_k^2-E\bigl(Y_k^2\bigr)
\bigr]E(K_{p,kk})+O_p\bigl(p^{r/2-1}\bigr) \bigr\}.
\end{equation}
Note that $c_p^2\leq l_p=O(p^r)$ under \textup{(C.2)}. Combined with
Proposition~\ref{i-cook}, we have
\begin{eqnarray*}
B_1-E(B_1)&=&O_p\bigl(n^{-7/2}
\bigr),
\\
B_3-E(B_3)&=&O_p\bigl(n^{-5/2}p^{r/2-1}\bigr),
\\
B_4-E(B_4)&=&O_p\bigl(p^{r/2-1}n^{-2}\bigr).
\end{eqnarray*}
Consequently, we have%
%
\begin{eqnarray}\label{B134}
&& \frac{n(n-1)^2}{(n-2)} \biggl\{\sum_{i=1,3}
\bigl[B_i-E(B_i)\bigr]-2\bigl(B_4-E(B_4)
\bigr)\biggr\}
\nonumber\\[-8pt]\\[-8pt]
&&\qquad =O_p\bigl(n^{-3/2}\bigr)+O_p
\bigl(p^{r/2-1}\bigr).\nonumber
\end{eqnarray}
Furthermore, by the results on $E(\mathcal{D}_k)$ in step~2 of the
proof of Proposition~\ref{i-cook}, we have
\[
E(\mathcal{D}_k)=\frac{1}{n(n-1)}E\bigl(Y_k^2
\bigr) E(K_{p,kk})+O\bigl(n^{-2}p^{-1}l_p^{1/2}
\bigr).
\]
Consequently, by $l_p=O(p^r)$, we have%
%
\begin{equation}
\label{EDk} \frac{n(n-1)^2}{(n-2)} E(\mathcal{D}_k)=\frac{n-1}{n-2}E
\bigl(Y_k^2\bigr)E(K_{p,kk})+O
\bigl(p^{r/2-1}\bigr).
\end{equation}
Combining (\ref{B2})--(\ref{EDk}), we have
\begin{eqnarray*}
&& \frac{n(n-1)^2}{(n-2)}\mathcal{D}_k
\\
&&\qquad =\frac{n(n-1)^2}{(n-2)} E(
\mathcal{D}_k)+ \frac
{n(n-1)^2}{(n-2)}\bigl[\mathcal{D}_k-E(
\mathcal{D}_k)\bigr]
\\
&&\qquad =\frac{n(n-1)^2}{(n-2)} E(\mathcal{D}_k)+ \frac{n(n-1)^2}{(n-2)}
\biggl(\sum_{i=1,2,3}\bigl[B_i-E(B_i)
\bigr]-2\bigl(B_4-E(B_4)\bigr) \biggr)
\\
&&\qquad = \frac{n-1}{n-2} E\bigl(Y^2_k
\bigr)E(K_{p,kk})+ \bigl\{Y_k^2-E
\bigl(Y_k^2\bigr) \bigr\} E(K_{p,kk})
\\
&&\quad\qquad{} +O_p\bigl(p^{r/2-1}\bigr)+O_p
\bigl(n^{-3/2}\bigr)
\\
&&\qquad = Y_k^2E(K_{p,kk})+\frac{1}{n-2} E
\bigl(Y_k^2\bigr)E(K_{p,kk})+O_p
\bigl(p^{r/2-1}\bigr)+O_p\bigl(n^{-3/2}\bigr)
\\
&&\qquad = Y_k^2+o_p(1),
\end{eqnarray*}
where the\vspace*{-1pt} last equation is from the fact that $E(X_{kj}^2)=1,
j=1,\ldots,p$ and $E(Y_k^2)=1$. Since $Y\sim N(0,1)$, we have
$\frac{n(n-1)^2}{(n-2)}\mathcal{D}_k\sim\chi^2(1)$; that is, $n^2
\mathcal{D}_k\sim\chi^2(1)$.
\end{pf*}
\end{appendix}

\section*{Acknowledgements}
We thank Professor Peter B\"uhlmann, Professor Peter Hall, an Associate
Editor and three referees for their constructive comments.

\begin{supplement}[id=suppA]
\stitle{Further proofs}
\slink[doi]{10.1214/13-AOS1165SUPP} 
\sdatatype{.pdf} \sfilename{aos1165\_supp.pdf}
\sdescription{The supplementary file contains the proofs of four
additional lemmas, Proposition~\ref{Dksample} and Theorem~\ref{Cook-u}.}
\end{supplement}

%

\printaddresses

\end{document}